\documentclass[reqno,a4paper,11pt]{amsart}
\usepackage[utf8]{inputenc}

\usepackage[english]{babel}

\usepackage[T1]{fontenc}
\usepackage[normalem]{ulem}

\usepackage{amssymb}
\usepackage{amsmath,amsfonts}
\usepackage[tikz]{bclogo}
\usetikzlibrary{arrows, arrows.meta}
\usepackage[top=3cm,bottom=3cm,left=2.5cm,right=2.5cm]{geometry}

\usepackage{lmodern}
\usepackage{textcomp}
\usepackage{mathtools, bm}
\usepackage[unq]{unique}
\usepackage{comment}
\usepackage[breaklinks]{hyperref}

\usepackage[numbers]{natbib}  

\usepackage{graphicx} 

\usepackage{amsthm}

\usepackage{subcaption}

\usepackage{cleveref}

\usepackage{bbm}

\usepackage{enumitem}
\usepackage{thmtools}
\usepackage{thm-restate}

\newtheorem{theorem}{Theorem}[section]
\newtheorem{lemma}[theorem]{Lemma}

\theoremstyle{definition}

\newtheorem{conjecture}[theorem]{Conjecture}

\DeclareMathOperator{\ind}{\mathbf{1}}

\newcommand{\EE}{\mathbb E}
\newcommand{\PP}{\mathbb P}

\newcommand{\dtv}{d_{TV}}

\title{Sprinkling for random regular graphs and Kim-Vu's sandwich conjecture}
\author{Julien Portier}
\address{Institute of Mathematics, EPFL, Lausanne, Switzerland}
\email{julienportier@laposte.net}
\date{}

\begin{document}

\begin{abstract}
    We show that the uniform random $d$-regular graph model on $n$ vertices is, in total variation distance, asymptotically close to the union of $d$ uniformly random edge-disjoint perfect matchings, provided that $n$ is even and $\omega(1) = d \leq n^{1/10}$. This confirms a conjecture by Hollom, Lichev, Mond, Portier, and Wang for this range of parameters.
    This result has a couple of interesting consequences.
    First, this shows that random $d$-regular graphs on $n$ vertices exhibit sprinkling properties for some range of parameters, thereby confirming a conjecture by Isaev, McKay, Southwell and Zhukovskii in that range.
    Second, we give a short alternative proof of the upper side of Kim-Vu's sandwich conjecture in the sparse regime.
\end{abstract}
\maketitle

\section{Introduction}

A recurring theme in random graph theory is to understand how different random graph models, or different instances of the same model, can be related to one another. One of the most useful examples is the technique of \emph{sprinkling}. If $G_p$ denotes the random subgraph of a fixed graph $G$ obtained by retaining each edge independently with probability $p$, then, whenever $1-p=(1-p_1)(1-p_2)$, one may generate $G_p$ as the union of two independent percolations with parameters $p_1$ and $p_2$. This simple observation makes it possible to reveal a random graph in several stages, and has proved extremely powerful in the fields of random graph theory and percolation. However, for random regular graphs, there is no analogous immediate such property, since the degree constraint introduces dependencies between the edges. \\

We write $G(n,d)$ for the uniform random $d$-regular graph on $[n]$, and denote its distribution by $\mu_d$. For two distributions on regular graphs, we write $\mu_{d_1}\oplus\mu_{d_2}$ for the distribution obtained by sampling $G_1\sim\mu_{d_1}$ and $G_2\sim\mu_{d_2}$ independently, conditioning on $G_1$ and $G_2$ being edge-disjoint, and taking their union. Motivated by the analogy with sprinkling, Isaev, McKay, Southwell and Zhukovskii~\cite{IMcKSZ23} conjectured that the distribution of $G(n,d_1)\oplus G(n,d_2-d_1)$ is essentially the one of $G(n,d_2)$ as long as $d_2$ tends to infinity.
\begin{conjecture}\label{conj:IMSZ-coupling-sum}
For $d_1,d_2\in [1,n-1]$ with $d_1\le d_2 =\omega(1)$, there exists a coupling of $G_{\oplus} \sim G(n,d_1) \oplus G(n,d_2-d_1)$ and $G\sim G(n,d_2)$ such that $\mathbb{P}(G_{\oplus} = G)=1-o(1)$.
\end{conjecture}
They verified this conjecture in several regimes, namely:
\begin{itemize}
\setlength\itemsep{0.3em}
    \item $d_1 = 1$, $d_2=\omega(1)$ and $d_2 = o(n^{1/3})$,
    \item $\min\{d_1,d_2-d_1,n-d_2\} = \omega(n/\log n)$,
    \item $d_1^2 (n/\min\{d_2-d_1,n-d_2+d_1\})^{2d_1}\le (\log n)/3$.
\end{itemize}
Moreover, Hollom, Lichev, Mond, Wang and the present author proved that, when $d_2 \leq n^{1/10}$, the distributions of $G(n,d_1)\oplus G(n,d_2-d_1)$ are asymptotically the same for all $d_1$ far from $1$ and $d_2$ (see Theorem 1.4 in \cite{hollom2025monotonicity}).
These results provide evidence that random regular graphs admit a genuine analogue of sprinkling despite the dependence between their edges. \\

In this paper we approach this question through a closely related model. Suppose that $n$ is even, and let $\nu_d=\mu_1 \oplus \dots \oplus \mu_1$ (with $d$ copies of $\mu_1$), i.e. the distribution obtained by taking the union of $d$ uniformly random edge-disjoint perfect matchings. 
Our first main result shows that, throughout a polynomial range of degrees, this model is close to the uniform random regular graph in total variation.

\begin{theorem}
\label{thm:Coupling-Decomposition-PerfectMatchings}
Let $n\to\infty$ through even integers, and let $d=d(n)$ satisfy $\omega(1) = d\leq n^{1/10}$.
Then
\[
 \dtv(\mu_d,\nu_d)=o(1).
\]
\end{theorem}

This proves, in the range $d\leq n^{1/10}$, Conjecture~1.6 of
Hollom, Lichev, Mond, Wang and the present author from
\cite{hollom2025monotonicity} in a strong form.
Previously, only a one-sided contiguity result was known in this range (see Theorem 1.5 in \cite{hollom2025monotonicity}).
Largely as a consequence of \Cref{thm:Coupling-Decomposition-PerfectMatchings}, we obtain the following sprinkling result, which confirms \Cref{conj:IMSZ-coupling-sum} of Isaev, McKay, Southwell and Zhukovskii in a new range of parameters.

\begin{theorem}
\label{thm:Coupling-Sprinkling}
Let $n\to\infty$ through even integers. Let
$1\leq d_1<d_2 \leq n^{1/10}$ be such that $d_2=\omega(1)$.
Then
\[
 \dtv\bigl(\mu_{d_1}\oplus\mu_{d_2-d_1},\mu_{d_2}\bigr)=o(1).
\]
Equivalently, there is a coupling of $G_{\oplus}\sim G(n,d_1)\oplus G(n,d_2-d_1)$ and $G\sim G(n,d_2)$ such that $\PP(G_{\oplus}=G)=1-o(1)$.
\end{theorem}

In this paper, we slightly abuse notation by using $G(n,\cdot)$ for three different random graph models: $G(n,d)$ denotes the uniform random $d$-regular graph on $[n]$, $G(n,p)$ the binomial Erd\H{o}s-R\'enyi random graph with edge probability $p$, and $G(n,M)$ the uniform random graph on $[n]$ with exactly $M$ edges. The intended meaning will always be clear from context. \\

\textbf{Kim-Vu's sandwich conjecture.} 
Due to the lack of independence between edges in the model of random regular graphs, properties which are relatively easy to study for the Erd\H{o}s-R\'enyi random model, such as Hamiltonicity, are much harder to obtain for the random regular graph model. 
This led Kim and Vu \cite{kim2004sandwiching} to formulate the "sandwich conjecture", which asserts that a random regular graph can be coupled between two Erd\H{o}s-R\'enyi random graphs of essentially the same density.
More precisely, they conjectured that, if $d=\omega(\log n)$, then there exist $p_-,p_+=(1+o(1))\frac{d}{n}$ and a coupling $(G_-,G,G_+)$ such that $G_-\sim G(n,p_-)$, $G\sim G(n,d)$, $G_+\sim G(n,p_+)$ and $\mathbb{P}\bigl(G_-\subseteq G\subseteq G_+\bigr)=1-o(1)$.
Besides being interesting in its own right, such a comparison allows one to transfer many properties from the Erd\H{o}s-R\'enyi model to random regular graphs.
Substantial progress towards this conjecture was made by Kim and Vu \cite{kim2004sandwiching}, Dudek, Frieze, Ruci{\'n}ski, and {\v{S}}ileikis \cite{dudek2017embedding}, Gao, Isaev and McKay \cite{gao2022sandwiching,gao2011kim}, Klimo{\v{s}}ov{\'a}, Reiher, Ruci{\'n}ski, and {\v{S}}ileikis \cite{klimovsova2023sandwiching}, before the conjecture was proved in full by Behague, Il'kovi{\v{c}}, and Montgomery \cite{behague2025proof}.
Using \Cref{thm:Coupling-Decomposition-PerfectMatchings}, we give an alternative short proof of the upper side of the sandwich in the sparse regime (which were the last remaining range and side before the work of \cite{behague2025proof}).

\begin{theorem}\label{thm:Kim-Vu}
Let $n\to\infty$ through even integers, and suppose
$\log n\ll d\le n^{1/10}$. Then there exists $p_+ =(1+o(1))\frac dn$ such that there is a coupling of $G_d\sim\mu_d$ and $G_+\sim G(n,p_+)$
satisfying
\[
 \mathbb P(G_d\subseteq G_+)=1-o(1).
\]
\end{theorem}

We remark that it would not be hard to remove the assumption that $n$ is even in the previous result by using the approach of Gao, Isaev, and Perez-Gimenez \cite{gao2025sandwiching} (see pages 28-29) or alternatively, Claim 4.26 from \cite{hollom2026approximate} (which can be proved for $d \leq n^{1/10}$ via the same proof contained in \cite{hollom2026approximate}). For simplicity, however, we restrict throughout to even $n$. \\

Our goal in giving this alternative proof of \Cref{thm:Kim-Vu} is two-fold. 
First we want to illustrate the usefulness of \Cref{thm:Coupling-Decomposition-PerfectMatchings} as a tool, as it implies the upper side of Kim-Vu's sandwich conjecture with relatively little additional work.
Moreover, as mentioned above, a one-sided contiguity version of
\Cref{thm:Coupling-Decomposition-PerfectMatchings} was previously proved by Hollom, Lichev, Mond, Wang and the present author in~\cite{hollom2025monotonicity}, and has already proved useful in \cite{hollom2026approximate} where the same authors proved a weak form of the Itai-Zehavi conjecture for random regular graphs.
Second, and maybe more importantly, we believe that our argument gives a natural explanation for why the upper side of Kim-Vu's sandwich conjecture holds: because random regular graphs are close to unions of random edge-disjoint perfect matchings, and those perfect matchings switch easily. \\

\textbf{Future directions of research.}
We believe that \Cref{thm:Coupling-Decomposition-PerfectMatchings} may have further applications, and its main current limitation is that it is proved only in the sparse regime.
The author conjectures that \Cref{thm:Coupling-Decomposition-PerfectMatchings} can be extended to $d=o(n)$, which would in particular imply Conjecture 1.6 of \cite{hollom2025monotonicity} in this range.
If this were true, we hope this would provide a useful tool for the study of random regular graphs.
For instance, if this conjecture turned out to be true, then our proof of \Cref{thm:Kim-Vu} can be adapted to a proof of the upper side to the whole range $d=o(n)$. We note, however, that the arguments from \cite{gao2025sandwiching} and \cite{hollom2026approximate} mentioned above to extend to $n$ odd would not immediately extend to this range.
The polynomial restriction on $d$ in \Cref{thm:Coupling-Decomposition-PerfectMatchings} comes from two sources. The first is Theorem 4.6 in \cite{mckay1985asymptotics}, which only provides a good estimate on the number of perfect matchings in the complement of a sparse graph $X$. The second is Corollary~3.2 of~\cite{hollom2025monotonicity} (see \eqref{eq:dtv-Oldpaper}), which shows that the edge-disjoint union of a random $d$-regular graph and a random perfect matching is roughly a random $(d+1)$-regular graph, and is currently only proved in this degree range. 
This second restriction is itself due to a result on the concentration of the number of perfect matchings in $G(n,d)$ (see Proposition 3.1 of \cite{hollom2025monotonicity}) which is only available in this range of $d$. If one could prove such a concentration for higher values of $d$, this would remove one of the main bottlenecks to extend \Cref{thm:Coupling-Decomposition-PerfectMatchings}. \\

\textbf{AI Declaration. }
The proofs of \Cref{thm:Coupling-Decomposition-PerfectMatchings} and \Cref{thm:Coupling-Sprinkling} are due to the author. 
However, the author first produced a completely different proof of \Cref{thm:Coupling-Decomposition-PerfectMatchings} which heavily relied on ChatGPT assistance. 
More precisely, the author used ChatGPT to prove that, for fixed $d$, we have $\dtv(\mu_d, \nu_d) = O(d^{-1/2})$, from which the author knew that this bound could be lifted to imply \Cref{thm:Coupling-Decomposition-PerfectMatchings} by using the "peeling off one matching at a time" result by Hollom, Lichev, Mond, Wang and the present author (see Corollary 3.2 in \cite{hollom2025monotonicity}, or \eqref{eq:dtv-Oldpaper}), which itself crucially builds upon the ideas of Gao \cite{Gao23}. 
The proof of $\dtv(\mu_d, \nu_d) = O(d^{-1/2})$ relied on a second moment computation, and the clever observation by ChatGPT that this second moment corresponds to counting regular graphs equipped with a colouring having some particular properties, which is exactly what the result of Borinsky, Meroni and Wiesmann~\cite{BMW2026} provides. However, its asymptotic formula depends on auxiliary analytic quantities, namely certain critical points and Hessian determinants, which are not particularly difficult to compute, although quite lengthy. 
While writing the details of those computations, the author realised that using the "retry coupling" idea from Lemma 5.3 in the work of Hollom, Lichev, Mond, Wang and the present author \cite{hollom2025monotonicity}, one could in fact get a completely different proof that is much simpler and shorter, and is the one presented here.
Once \Cref{thm:Coupling-Decomposition-PerfectMatchings} was established, the author believed that one could obtain a simpler proof of the upper side of Kim-Vu's sandwich conjecture using \Cref{thm:Coupling-Decomposition-PerfectMatchings} as a black box.
The main motivation for this was that the proof of Behague, Il'kovi{\v{c}}, and Montgomery \cite{behague2025proof} relies on switching arguments, and therefore working with perfect matchings could make those switching arguments much simpler. 
The author then asked ChatGPT to pursue this idea. Through several iterations, the model produced a rather complicated proof that the author did not check in detail. At this stage, no mathematical input from the author was needed beyond prompting the model to continue proving the technical lemmas it identified as necessary.
Upon inspection, this new proof relied on complicated spectral arguments to control the in- and out-degrees of the switching graph. The author believed those could be replaced by a more direct combinatorial and probabilistic approach, while keeping the same switching strategy.
When prompted to do so, the model was able to produce such a proof, which is the one presented here, up to a few simplifications and revisions by the author.
The author independently verified and revised all arguments and computations appearing in the final manuscript, and takes full responsibility for the content of the paper. \\

\textbf{Acknowledgements.} The author would like to thank Lawrence Hollom, Lyuben Lichev, Adva Mond, and Yiting Wang for valuable discussions identifying that proving the result $\limsup_{d \to \infty}\limsup_{n \to \infty} \dtv(\mu_d, \nu_d) =0$ was enough to show \Cref{thm:Coupling-Decomposition-PerfectMatchings}, and attempts toward proving this result.

\section{Outline of proofs}

We first give some ideas for the proof of \Cref{thm:Coupling-Decomposition-PerfectMatchings}, as \Cref{thm:Coupling-Sprinkling} will mostly follow from \Cref{thm:Coupling-Decomposition-PerfectMatchings} and the methods used to prove it. 
First, it has been proved by Hollom, Lichev, Mond, Wang, and the author \cite{hollom2025monotonicity}, building on ideas of Gao \cite{Gao23}, that in our degree range,
\begin{align}
\label{eq:dtv-Oldpaper}
    \dtv(\mu_{d+1},\mu_d\oplus\mu_1)=O(d^{-1.1}),
\end{align}
where we recall that $\mu_r$ is the uniform distribution on $r$-regular graphs (and so $\mu_1$ is the uniform distribution on perfect matchings) and
$\oplus$ denotes independent union conditioned on edge-disjointness.
Since the errors in \eqref{eq:dtv-Oldpaper} are summable, to establish our desired result, we may essentially peel off one matching at a time. 
Thus, it is sufficient to prove a good bound for small $d$, which we can then lift to growing degrees using the above result. 
In fact, still by the argument of peeling off one matching at a time, it is sufficient to obtain that, for fixed $m \geq 3$,
\begin{align*}
    \lim_{k \to \infty} \limsup_{n \to \infty} \dtv(\mu_m \oplus \nu_{km-m}, \nu_{km }) =0.
\end{align*}
To prove this, we adapt the sequential coupling idea of  Lemma 5.3 of \cite{hollom2025monotonicity}. The main idea is to write $\nu_{km}$ as $\nu_m \oplus \dots \oplus \nu_m$ where $ \nu_m$ appears $k$ times, and successively attempt to couple $\mu_m$ with each of the $k$ blocks, until a coupling succeeds.
We then move the matching block to the front and use the
remaining blocks to construct the complementary union.
The main point is to control the bias introduced by
unsuccessful attempts, so that the resulting distribution
is close to $\mu_m\oplus\nu_{km-m}$.
The full argument is written in \Cref{lem:sequential-coupling}. \\

We now give some ideas for the proof of \Cref{thm:Kim-Vu}. 
Our starting point is the coupling approach of Gao, Isaev, and McKay \cite{gao2022sandwiching}, which is also used in the final proof of Kim-Vu's conjecture by Behague, Il'kovi{\v{c}}, and Montgomery \cite{behague2025proof}.
In this approach, we work with a random graph $F_r$ obtained from $G\sim\mu_d$ by adding $r$ uniformly random edges of its complement. 
We start from $F_{\binom n2-m}=K_n$, where $m=dn/2$, and construct
a deletion process that goes from $F_r$ to $F_{r-1}$, and the goal is to show that this process can be coupled closely with the ordinary uniform edge deletion.
Note that the deletion probability of an edge $e$ in $F_r$ is proportional to
the number of $d$-regular subgraphs of $F_r$ avoiding $e$.
Thus, to show that the two deletion rules are close, the key task is to show
that every edge belongs to roughly the same proportion of
$d$-regular subgraphs.
This is the key technical challenge that \cite{gao2022sandwiching} and \cite{behague2025proof} encounter in their proofs. \\

However, by \Cref{thm:Coupling-Decomposition-PerfectMatchings}, we may replace $\mu_d$ by $\nu_d$ (which we recall is the measure obtained by taking the union of $d$ random edge-disjoint perfect matchings).
We therefore need to compare the numbers of $d$ edge-disjoint perfect matchings subgraphs containing
and avoiding a fixed edge $e$, rather than the numbers of $d$-regular
subgraphs. This makes the switching argument substantially simpler, as one can fix one of the perfect matchings $M_i$ and consider a cycle whose edges alternate between edges of $M_i$ and edges unused by the entire tuple. Flipping along such a cycle, that is, removing from $M_i$ the matching edges of the cycle and replacing them by the unused ones, produces another perfect matching while leaving all the other matchings unchanged.
To count those switches, we let $\rho=r/(\binom n2-m)$ be the density of the
added edges, and fix some integer $k$ and an edge $e$.
If $e$ belongs to the tuple, then its colour is fixed. An alternating cycle of length $2k$ containing $e$ has roughly $n^{k-1}$ possible choices: the matching steps are forced, while at each unused step there are roughly $n$ choices, except for the last one, for which we have only one choice. Since all $k$ unused edges must be present, this gives roughly $2n^{k-1}\rho^k$ switchings removing $e$, where the factor $2$ takes into account the two possible orientations of the cycle. Similarly, if $e$ does not belong to the tuple, we may insert it into any of the $d$ matchings, and $e$ itself provides one of the unused edges of the cycle, giving roughly $2dn^{k-2}\rho^{k-1}$ switchings inserting $e$ (see \Cref{lem:pseudorandom-properties} for the formal proof of both estimates).
By double-counting one finds that the
proportion $x_e(F)$ of tuples containing $e$ satisfies
\[
 x_e(F)\approx
 \frac{2dn^{k-2}\rho^{k-1}}
      {2n^{k-1}\rho^k+2dn^{k-2}\rho^{k-1}}
 =\frac{d}{d+n\rho}
 \approx\frac{m}{m+r},
\]
and thus the deletion probabilities, which are proportional to $1-x_e(F)$,
are nearly uniform.
This is the main idea of the proof. A few additional technicalities arise in making the argument rigorous. First the above switching estimates hold only with high probability, and thus one must take into account those exceptional events and show that they do not contribute many switchings, which can be done by proving upper bounds on the degrees in $F_r$. Second the error between the deletion probabilities is too large to couple the two deletion processes so that they delete the same edge at every step. Instead, we delete the same edge at most steps of the coupled process, and show that the exceptional steps are sufficiently few to be absorbed by a small surplus of edges in the uniform process.

\section{Sprinkling for random regular graphs}

Recall that for probability measures on graphs with vertex set $[n]$, we write
$P_1\oplus\cdots\oplus P_k$ for the law of their independently
sampled union, jointly conditioned on the graphs being pairwise
edge-disjoint. 
For convenience, we write $Q^{\oplus k}$ for $k$ copies of $Q$.

The following result is our main tool to prove \Cref{thm:Coupling-Decomposition-PerfectMatchings}. Both its statement and proof are largely inspired by Lemma~5.3 of \cite{hollom2025monotonicity}, although we adopt a slightly different formulation and presentation, which we find more natural and easier to understand.

\begin{lemma}\label{lem:sequential-coupling}
Fix $r\ge1$, and let $\alpha_n,\beta_n$ be relabelling-invariant
probability measures on the simple $r$-regular graphs on $[n]$, where $n$ is even.
If $\alpha_n$ is contiguous with respect to $\beta_n$, then
\[
 \lim_{k\to\infty}\limsup_{n\to\infty}
 \dtv\bigl(\alpha_n\oplus\beta_n^{\oplus(k-1)},
                  \beta_n^{\oplus k}\bigr)=0.
\]
\end{lemma}

\begin{proof}

Let $P_n,Q_n$ be the distributions induced by $\alpha_n,\beta_n$,
respectively, on the isomorphism classes of $r$-regular graphs.
We first reduce the problem to a coupling of these classes.
For fixed $k$, let $K_{n,k}$ be the kernel that, given $k$ classes,
chooses uniform labelled representatives conditioned on being
pairwise edge-disjoint, and outputs their union. By the argument of
\cite[Lemma~5.2]{hollom2025monotonicity}, using the uniform overlay
estimate \cite[Lemma~2.10]{hollom2025monotonicity}, we have that
\begin{align*}
 \dtv\bigl(K_{n,k}(P_n\otimes Q_n^{\otimes(k-1)}),
       \alpha_n\oplus\beta_n^{\oplus(k-1)}\bigr)=o(1), \qquad \text{and} \qquad \dtv\bigl(K_{n,k}(Q_n^{\otimes k}),\beta_n^{\oplus k}\bigr)=o(1),
\end{align*}
where $S\otimes T$ denotes the joint law of two independent
samples with respective distributions $S$ and $T$, and
$S^{\otimes k}$ denotes the law of $k$ independent samples from $S$.
Thus by triangle inequality we have
\begin{align}
\label{eq:asp-mixed}
    \dtv\bigl(\alpha_n\oplus\beta_n^{\oplus(k-1)},
                  \beta_n^{\oplus k}\bigr) \leq \dtv\bigl(K_{n,k}(P_n\otimes Q_n^{\otimes(k-1)}),K_{n,k}(Q_n^{\otimes k})\bigr) +o(1).
\end{align}
Consequently, it suffices to show that, for every $\varepsilon>0$,
all sufficiently large $k$ and then all sufficiently large $n$ admit
a coupling of $Y\sim Q_n^{\otimes k}$ with a random permutation
$Z$ of $Y$, whose law $\Pi$ satisfies
\begin{equation}\label{eq:target}
 \dtv\bigl(\Pi,P_n\otimes Q_n^{\otimes(k-1)}\bigr)<\varepsilon.
\end{equation}
Indeed, since $K_{n,k}(\Pi)=K_{n,k}(Q_n^{\otimes k})$, the claimed bound
on the graph unions follows by \eqref{eq:asp-mixed}. 

We now further reduce the problem by showing that we may assume $P_n\le M Q_n$ for a fixed $M\ge1$.
Fix $\eta\in(0,1/2)$. By contiguity, there exists a sufficiently large constant
$L$ such that, for all sufficiently large $n$,
\[
 P_n(B_n)\ge 1-\eta,
\]
where $B_n=\{G:P_n(G)\le LQ_n(G)\}$.
Indeed, otherwise, since $Q_n(B_n^c)\le1/L$, by taking $L\to\infty$
and a corresponding subsequence of $n$ would contradict contiguity.
The measure $\widetilde P_n=P_n(\,\cdot\mid B_n)$ satisfies
\[
 \dtv(P_n,\widetilde P_n)\le\eta,
 \qquad \widetilde P_n\le M Q_n,
\]
where  $M=2L$.
Thus replacing $P_n$ by $\widetilde P_n$ costs at most $\eta$
in \eqref{eq:target}. Since $\eta$ is arbitrary, it is therefore enough to
prove \eqref{eq:target} under the assumption that $P_n\le M Q_n$, with $M$ independent of $n$.

Choose an integer $m\ge2M$ and, for convenience, suppress $n$ from the notation. Let $a=1-\frac1m$ and let $R$ be the probability measure such that
\begin{align*}
    Q=m^{-1}P+aR.
\end{align*}
We now construct the coupling, following the idea of Lemma 5.3 from \cite{hollom2025monotonicity}.
Sample $X\sim P$. At each step, until the first accepted match,
accept with probability $1/m$, independently of the previous
choices. On acceptance, set $Y_i=X$, and on rejection, sample
$Y_i\sim R$ independently of $X$ and all previous rejected
samples. 
After acceptance, sample the remaining blocks
independently from $Q$. 

Since $m^{-1}P+aR=Q$, we have that
$Y=(Y_1,\ldots,Y_k)\sim Q^{\otimes k}$.
Move the first accepted block to the front to obtain $Z$.
If all attempts fail, set $Z=Y$. Thus, $Z$ is always
obtained from $Y$ by a move-to-front permutation.
The law $\Pi$ of $Z$ is
\begin{equation}
\label{eq:retry-law}
 \Pi=a^kR^{\otimes k}
 +\sum_{j=1}^k\frac{a^{j-1}}m
       P\otimes R^{\otimes(j-1)}\otimes Q^{\otimes(k-j)} = S^{\otimes k}+\frac1m P\otimes
       \sum_{j=0}^{k-1}S^{\otimes j}\otimes Q^{\otimes(k-1-j)},
\end{equation}
where $S=aR=Q-\frac1m P$.
Our goal is now to compare this distribution with $A=P\otimes Q^{\otimes(k-1)}$, and show that they are close in total variation distance.
Since $Q=S+P/m$, a simple expansion gives
\[
A=P\otimes Q^{\otimes(k-1)}
=P\otimes S^{\otimes(k-1)}
+\frac1m\sum_{j=0}^{k-2}
P\otimes S^{\otimes j}\otimes P
\otimes Q^{\otimes(k-2-j)}.
\]
Subtracting this from the expression for $\Pi$ in \eqref{eq:retry-law}, and using again
$S+P/m=Q$, we obtain
\begin{equation}\label{eq:telescoping}
 \Pi-A=(Q-P)\otimes S^{\otimes(k-1)} +\frac1m P\otimes\sum_{j=0}^{k-2}
    S^{\otimes j}\otimes(Q-P)\otimes Q^{\otimes(k-2-j)}.
\end{equation}

For every graph $G$ in the support of $Q$, we let $f(G)=P(G)/Q(G)$ and $\tau=\EE_Q f^2$. Since $P\le MQ$,
we have $\EE_Q f=1$ and $1\le\tau\le M$.
Let $G_1,\ldots,G_{k-1}$ be independent with law $Q$ and for $ 1\le j\le k-1$, we let
\[
 D_j=\bigl(1-f(G_j)\bigr)
       \prod_{i=1}^{j-1}\left(1-\frac{f(G_i)}m\right).
\]
Summing \eqref{eq:telescoping} in absolute value over each graph $G$, and using that the total mass of $S$ is $a$, we get
\begin{equation}\label{eq:tv-telescope}
 2\dtv(\Pi,A)
 \le
 2a^{k-1}
 +\frac1m\EE\left|\sum_{j=1}^{k-1}D_j\right|.
\end{equation}
Note that $\EE[D_iD_j]=0$ for $i\ne j$. Moreover, we have
\begin{align*}
    \EE D_j^2=(\tau-1)\left(\EE_Q\left(1-\frac fm\right)^2\right)^{j-1} = (\tau-1)\left(1-\frac2m+\frac{\tau}{m^2}\right)^{j-1} \leq (\tau-1)a^{j-1}.
\end{align*}
Hence, by Cauchy-Schwarz we have
\[
 \frac1m\EE\left|\sum_{j=1}^{k-1}D_j\right|
 \le\frac1m\sqrt{\sum_{j=1}^{k-1}\EE D_j^2}
 \le\frac1m\sqrt{(\tau-1)\sum_{j=0}^\infty a^j}
 \le\sqrt{\frac{M-1}{m}}.
\]
Plugging the previous inequality into \eqref{eq:tv-telescope}, we deduce that
\begin{equation}\label{eq:retry-bound}
 \dtv\bigl(\Pi,P\otimes Q^{\otimes(k-1)}\bigr)
 \le \left(1-\frac1m\right)^{k-1}
      +\frac12\sqrt{\frac{M-1}{m}}.
\end{equation}
By letting $k\to\infty$, and then letting
$m\to\infty$, this shows that \eqref{eq:target} is satisfied, and thus finishes the proof.
\end{proof}

\Cref{thm:Coupling-Decomposition-PerfectMatchings} now essentially follows from \Cref{lem:sequential-coupling} by applications of Corollary 3.2 of \cite{hollom2025monotonicity}.

\begin{proof}[Proof of \Cref{thm:Coupling-Decomposition-PerfectMatchings}.]
    Fix $m,k \geq 3$, and let $d \geq mk$. For $m\leq i\leq d$, define
\[
 \lambda_i=\mu_i\oplus\nu_{d-i}.
\]
Note that by triangle inequality, we have 
\begin{align}
\label{eq:Triangle-Ineq-Tech}
    \dtv(\mu_d,\nu_d) \leq \dtv(\mu_d,\lambda_m) + \dtv(\lambda_m,\nu_d).
\end{align}
For every $i\in\{m,\ldots,d-1\}$, we have $\lambda_{i+1}
 = \mu_{i+1}\oplus\nu_{d-i-1}$ and by associativity $\lambda_i
 = (\mu_i\oplus\mu_1)\oplus\nu_{d-i-1}$.
Thus combining Corollary~2.9(a) and Corollary~3.2 from \cite{hollom2025monotonicity}, we obtain
\[
 \dtv(\lambda_i,\lambda_{i+1})
 \leq C\,\dtv(\mu_i\oplus\mu_1,\mu_{i+1})
 \leq C' i^{-1.1},
\]
for absolute constants $C,C'>0$. Consequently, by triangle inequality
\begin{equation}
\label{eq:fixed-to-growing-second}
 \dtv(\lambda_m,\mu_d)
 \leq \sum_{i=m}^{d-1}\dtv(\lambda_i,\lambda_{i+1})
 \leq C'\sum_{i=m}^{\infty}i^{-1.1}
 =O(m^{-0.1}).
\end{equation}
By Corollary~2.9(a) from \cite{hollom2025monotonicity}, we also find that
\begin{align*}
    \dtv(\lambda_m,\nu_d) \leq C \dtv(\mu_m \oplus \nu_{(k-1)m}, \nu_{km}).
\end{align*}
By writing $\nu_{(k-1)m}= \nu_m \oplus \dots \oplus \nu_m$ with $k-1$ copies of $\nu_m$, and similarly for $\nu_{km}$, we may apply \Cref{lem:sequential-coupling} (as it is well known from \cite{janson1995random} and \cite{MRRW97} that $\mu_m$ and $\nu_m$ are contiguous for fixed $m$) and get that
\begin{align}
\label{eq:second-term-decomp}
    \limsup_{n \to \infty} \dtv(\lambda_m,\nu_d) \leq C \limsup_{n \to \infty} \dtv(\mu_m \oplus \nu_{(k-1)m}, \nu_{km}) \leq \delta_{k,m},
\end{align}
where $\lim_{k\to\infty} \delta_{k,m} =0$.
Thus, by plugging in \eqref{eq:fixed-to-growing-second} and \eqref{eq:second-term-decomp} into \eqref{eq:Triangle-Ineq-Tech}, we find that
\begin{align*}
    \limsup_{n \to \infty} \dtv(\mu_d,\nu_d) \leq O(m^{-0.1})+ \delta_{k,m}.
\end{align*}
Taking $k\to\infty$, and then $m \to\infty$ gives the desired conclusion.
\end{proof}

We now prove \Cref{thm:Coupling-Sprinkling}. Under the assumption $\min\{d_1,d_2-d_1\}=\omega(1)$, the result follows easily from \Cref{thm:Coupling-Decomposition-PerfectMatchings}. However, to cover the full range of $d_1$, we make use of \Cref{lem:sequential-coupling}, following an argument similar to that in the proof of \Cref{thm:Coupling-Decomposition-PerfectMatchings}.

\begin{proof}
Let $b=d_2-d_1$. By symmetry we may assume that $d_1\le b$,
and thus that $b=\omega(1)$. 
Fix an integer $m\ge3$.
For every fixed $1\le s\le m$, the measures
$\mu_s\oplus\nu_2$ and $\nu_{s+2}$ are contiguous (again by \cite{janson1995random} and \cite{MRRW97}).
Applying \Cref{lem:sequential-coupling} to $\alpha_n=\mu_s\oplus\nu_2 $ and $\beta_n=\nu_{s+2}$, 
and adding perfect matchings using
Corollary~2.9(a) of \cite{hollom2025monotonicity}, we obtain
\begin{align}
\label{eq:maxs}
    \max_{1\le s\le m}
 \dtv(\mu_s\oplus\nu_{d_2-s},\nu_{d_2})=o(1).
\end{align}
Let $r=\min\{d_1,m\}$. By triangle inequality and the same
telescoping argument as in the proof of
\Cref{thm:Coupling-Decomposition-PerfectMatchings},
using Corollaries~2.9(a) and~3.2 of
\cite{hollom2025monotonicity}, we have
\begin{align*}
\dtv(\mu_{d_1}\oplus\mu_b,\mu_{d_2})
&\le C\dtv(\mu_b,\nu_b)
   +C\sum_{j=r}^{d_1-1}j^{-1.1} +\dtv(\mu_r\oplus\nu_{d_2-r},\nu_{d_2})
   +\dtv(\nu_{d_2},\mu_{d_2}).
\end{align*}
The first and fourth terms are $o(1)$ by \Cref{thm:Coupling-Decomposition-PerfectMatchings}.
The third term is also $o(1)$ by \eqref{eq:maxs}.
The second term is at most $C\sum_{j=m}^{\infty}j^{-1.1} =O(m^{-0.1})$.
Thus
\begin{align*}
\dtv(\mu_{d_1}\oplus\mu_b,\mu_{d_2}) = O(m^{-0.1})+o(1),
\end{align*}
and by letting $m \to \infty$, we obtain the desired result.
\end{proof}

\section{An alternative proof of the upper side of Kim-Vu's sandwich conjecture for sparse random graphs}

We start by introducing some notation. 
For an ordered tuple $\mathcal T=(M_1,\ldots,M_d)$ of edge-disjoint perfect matchings on $[n]$,
we write $K(\mathcal T)$ for the graph obtained by
taking their union, i.e. $E(K(\mathcal T))=\bigcup_{i=1}^d E(M_i)$.
We let $\mathcal M$ be uniform among all ordered $d$-tuples of
edge-disjoint perfect matchings of $K_n$, and then choose $A$ uniformly
among the $r$-subsets of $E(K_n)\setminus E(K(\mathcal M))$.
Let $F$ be the graph with edge set $E(K(\mathcal M)) \cup A$, and let $\pi_r$ be the distribution of $F$. 
Set $m=dn/2$ and $a=2r/n$.
For $i\in[d]$, $x,y\in[n]$ and $t\ge1$, let $X_t^{(i)}(x,y;\mathcal M,F)$ be the number of alternating paths
$x=v_0,v_1,\ldots,v_{2t}=y$ such that $v_0,\ldots,v_{2t-1}$ are
distinct, $v_{2t}\notin\{v_1,\ldots,v_{2t-1}\}$, $v_{2j}v_{2j+1}\in M_i$ and $v_{2j+1}v_{2j+2}\in E(F)\setminus E(K(\mathcal M))$ for every $0\le j<t$.
We now state our first lemma.

\begin{lemma}
\label{lem:pseudorandom-properties}
Let $\mathcal{M}$ and $F$ be defined as above.
Suppose $d\le n^{1/10}$ and $a\gg\log n$, and let $k=20 \log n$. 
Then there exists $\eta=O(\sqrt{\frac{\log n}{a}})$ such that with
probability at least $1-n^{-200}$, simultaneously for every colour
$i$, every $t\in\{k-1,k\}$ and every pair of vertices $x,y$ with $\{x,y\}\notin E(M_i)$, we have
\begin{align*}
 X_t^{(i)}(x,y;\mathcal M,F)
   =(1\pm\eta)\frac{a^t}{n},
\end{align*}
and for every vertex $v$,
\begin{align*}
 |d_F(v)-(d+a)|\le\eta a.
\end{align*}
\end{lemma}

\begin{proof}
We fix $\mathcal{M}$ for the rest of the proof.
We start by proving the first inequality. 
For simplicity, we prove it for $t=k$, the case $t=k-1$ being similar.
Let $\eta=C \sqrt{\frac{\log n}{a}}$ for some large constant $C$.
Let $A_\rho$ be the edge set obtained by taking each edge in $E(K_n)\setminus E(K(\mathcal M))$ independently with probability $\rho=r/L$, where $L= \binom{n}{2}-m$. Note that $A_\rho$ conditional on $|A_\rho|=r$ is exactly a uniform $r$-subset
of  $E(K_n)\setminus E(K(\mathcal M))$. Since $|A_\rho|$ follows a $\operatorname{Bin}(L,r/L)$ distribution, we have
\[
 \PP(|A_\rho|=r)\ge\frac1{L+1}\ge n^{-2}.
\]
Thus we may simply prove the desired inequalities for $A_\rho$ instead of $A$, with failure probability $n^{-202}$ instead of $n^{-200}$. \\

Fix a colour $i\in[d]$, and vertices $x,y$ such that $\{x,y\}\notin E(M_i)$.
Let $X=X_k^{(i)}(x,y;\mathcal M,K(\mathcal M)\cup A_\rho)$
and $\mu=\EE [X]$, where the expectation is over the randomness of $A_\rho$.
We start by estimating $\mu$. 
Since the matching steps are forced,
there are $k-1$ vertex choices, each with $n-O(d+k)$
possibilities, as we cannot repeat vertices, and must avoid edges of $K(\mathcal M)$. 
Thus the number of possible paths is
\[
n^{k-1}\left(1+O\left(\frac{k(d+k)}{n}\right)\right).
\]
Since each path requires $k$ distinct edges of $A_\rho$, and
$a=(n-1-d)\rho$ by definition, we obtain
\begin{align}
\label{eq:Value-mu}
    \mu=n^{k-1}\rho^k
\left(1+O\left(\frac{k(d+k)}{n}\right)\right)
=\frac{a^k}{n}(1+o(\eta)),
\end{align}
where we used $k=O(\log n)$ and $d\le n^{1/10}$. \\

We treat separately the lower and upper tails, while the computations share a similar idea.
We start with the lower tail. We will apply Janson's inequality, and thus we focus on estimating
\begin{align*}
    D= \sum_{\gamma, \gamma': \gamma \cap \gamma' \neq \emptyset} \mathbb{E}[\ind_{\gamma \subset A_\rho} \ind_{\gamma' \subset A_\rho}],
\end{align*}
where the sum is over alternating paths $\gamma,\gamma'$ of length
$2k$ from $x$ to $y$, whose edges alternate between $M_i$ and
$E(K_n)\setminus E(K(\mathcal M))$, starting with an edge of $M_i$, and whose vertices, except possibly the endpoints, do not repeat.
We slightly abuse notation by identifying each path with its set
of edges outside $K(\mathcal M)$ when writing
$\gamma\subset A_\rho$ or $\gamma\cap\gamma'$.
We estimate $D$ via a simple combinatorial argument.
First, note that 
\begin{align*}
    D= \mu + D_{proper},
\end{align*}
where $D_{proper}= \sum_{\gamma \neq \gamma': \gamma \cap \gamma' \neq \emptyset} \mathbb{E}[\ind_{\gamma \subset A_\rho} \ind_{\gamma' \subset A_\rho}]$.
Since we have already computed $\mu$, we focus on estimating $D_{proper}$.
For $\gamma$, we have at most $n^{k-1}$ possibilities. Suppose $\gamma$ and $\gamma'$ agree on exactly $\ell \neq k$ edges, and that those $\ell$ edges consist of $c$ internal components of $\gamma'$, where a component is a maximal sequence of
consecutive non-matching steps of $\gamma'$ shared
with $\gamma$, and a component is internal if it does not contain the first nor the last non-matching step of $\gamma'$.
Then we have $\binom{k}{c}$ choices to decide where those components start on $\gamma'$, and at most $\binom{\ell+1}{c+1}$ choices for the lengths of those components (including the ones containing potentially the first or last non-matching step). 
For each vertex on the path $\gamma'$, we have at most $n$ possibilities, except when we start an internal component, where we have $2k$ choices, and inside an internal component where we have to follow $\gamma$, and thus only $1$ possibility.
Therefore
\begin{align*}
    D_{proper} &\leq  \sum_{\ell = 1}^{k} \sum_{c=0}^{\ell} n^{k-1}\binom{k}{c}\binom{\ell+1}{c+1} (2k)^c n^{k-1-c-\ell} \rho^{2k-\ell} \\
    &\leq  (n^{k-1}\rho^{k})^2 \sum_{\ell = 1}^{k}  2 \left(\frac{2}{n\rho}\right)^{\ell}\sum_{c=0}^{\ell}\binom{k}{c} \left(\frac{2k}{n}\right)^c. 
\end{align*}
Since $\sum_{c=0}^{\ell}\binom{k}{c} \left(\frac{2k}{n}\right)^c \leq (1+\frac{2k}{n})^k \leq e^{2k^2/n}=O(1)$ and $\sum_{\ell = 1}^{k}  \left(\frac{2}{n\rho}\right)^{\ell} =O(\frac{1}{n\rho})$, we get
\begin{align*}
    D_{proper}= O \left((n^{k-1}\rho^{k})^2 \frac{1}{n\rho}\right).
\end{align*}
Since $\rho = r/L \geq a/n$ and $k \geq 10 \log n$, we also have that $\mu=O \left((n^{k-1}\rho^{k})^2 \frac{1}{n\rho}\right)$, and thus by Janson's inequality (see \cite{janson1990poisson}) we get
\begin{align*}
    \mathbb{P}(X \leq \mu (1- \eta) ) \leq \exp\left(- \frac{\eta^2 \mu^2}{2D}\right) \leq  \exp\left(- \Omega(\eta^2 \rho n )\right) \leq n^{-210}.
\end{align*}
Plugging in the value of $\mu$ from \eqref{eq:Value-mu}, and then taking a union bound over all colours, lengths and endpoint pairs, we obtain the desired result. \\

A similar computation treats the case of the upper tail. We will perform a high moment computation.
More precisely, for $j \leq n\rho/100$, we show an upper bound on 
\begin{align*}
    \mathbb{E}[X^j]=\sum_{\gamma_1, \dots, \gamma_j} \mathbb{E}[\ind_{\gamma_1 \subset A_\rho} \dots \ind_{\gamma_j \subset A_\rho}].
\end{align*}
We prove the bound inductively. Suppose that $\gamma_1,\ldots,\gamma_{j-1}$ have already been fixed, and let $U=\gamma_1\cup\cdots\cup\gamma_{j-1}$. 
To count the possible choices for $\gamma_j$, we first choose the positions at which it uses edges of $U$. Suppose that the overlap has size $1 \leq \ell<k$ and consists of $c$ internal components. 
As in the lower-tail argument, we have at most $\binom{\ell+1}{c+1}\binom{k}{c}$ choices for the overlap pattern.
For each internal component, the first overlapping edge can be chosen in at most $2k(j-1)$ ways, since $U$ contains at most $k(j-1)$ edges. Once inside an overlap component, each other overlapping edge has at most $j-1$ choices, since $U$ is the union of $j-1$ matchings.
Therefore, we obtain
\begin{align*}
    \mathbb{E}[X^j] &\leq  \mathbb{E}[X^{j-1} ] ( \mu+(2(j-1))^{k} + \sum_{\ell = 1}^{k} \sum_{c=0}^{\ell}\binom{\ell+1}{c+1}\binom{k}{c} (2k)^c (2j)^{\ell} n^{k-1-c-\ell} \rho^{k-\ell} ),
\end{align*}
where $\mu$ and $(2(j-1))^{k}$ are respectively the terms corresponding to an overlap of $\ell=0 $ and $\ell=k$ edges.
Now 
\begin{align*}
    \sum_{\ell = 1}^{k} \sum_{c=0}^{\ell}\binom{\ell+1}{c+1}\binom{k}{c} (2k)^c (2j)^{\ell} n^{k-1-c-\ell} \rho^{k-\ell} &\leq  2n^{k-1} \rho^k \sum_{\ell = 1}^{k}  \left(\frac{4j}{n\rho}\right)^{\ell} \sum_{c=0}^{\ell}\binom{k}{c}\left(\frac{2k}{n}\right)^c \\
    &\leq 2n^{k-1} \rho^k \sum_{\ell = 1}^{k} \left(\frac{4j}{n\rho}\right)^{\ell} \left(1+\frac{2k}{n}\right)^k \\
    &\leq  100\mu \frac{j}{n\rho},
\end{align*}
where we have used that $(1+\frac{2k}{n})^k \leq e^{2k^2/n} \leq e^2$, and that $j \leq n\rho/100$.
Thus it follows that, as long as $j\leq n\rho/100$, we have
\begin{align*}
    \mathbb{E}[X^j] &\leq \mu^j \prod_{i=1}^{j-1} \left(1+\frac{100(i+1)}{n\rho}+\frac{(2i)^k}{\mu}\right) \leq \mu^j  \exp \left(\frac{100j^2}{n\rho} + \frac{2^k j^{k+1}}{\mu} \right) \leq 2 \mu^j  \exp \left(\frac{100j^2}{n\rho}  \right),
\end{align*}
where we have used that $\frac{2^k j^{k+1}}{\mu}=o(1)$, since $k = 20 \log n$.
Therefore, by Markov's inequality, for $q= \eta n\rho/1000 $, we obtain
\begin{align*}
    \mathbb{P}\bigl(X\geq (1+\eta)\mu \bigr) \leq
    \frac{\mathbb{E}[X^q]}
    {(1+\eta)^q\mu ^q} \leq
    2 \exp\left(
        100\frac{q^2}{n\rho}
        -q\log(1+\eta)
    \right) \leq 2 \exp \left(
        -\Omega(\eta^2 n\rho)
    \right) \leq n^{-210}.
\end{align*}
Once again, plugging in the value of $\mu$ from \eqref{eq:Value-mu}, and then taking a union bound over all colours, lengths and endpoint pairs, we obtain the desired result. \\

The second inequality can be proved via a simple application of Chernoff's inequality.
Indeed, each degree of $A_\rho$ is binomial with mean $a$. Thus, for $0<\eta<1$,
Chernoff's inequality and a union bound give
\[
 \PP\left(\max_v|d_{A_\rho}(v)-a|>\eta a\right)
 \le2n\exp(-\eta^2a/3) \leq n^{-210}.
\]
This finishes the proof.
\end{proof}

Given a graph $F$, we let $\mathcal F_d(F)$ be the set of ordered $d$-tuples of edge-disjoint perfect matchings contained in $F$.
For an edge $e$, we let
\[
 x_e(F)=
 \frac{|\{\mathcal T\in\mathcal F_d(F):e\in K(\mathcal T)\}|}
      {|\mathcal F_d(F)|}.
\]
Next we show that, if $F$ follows a $\pi_r$ distribution, then typically all the $x_e(F)$ for $e \in E(F)$ have almost same value.
\begin{lemma}\label{lem:edge-marginals}
Let $n$ be even, let $1\le d\le n^{1/10}$, and let $r$ be an integer with
$0\le r\le \binom{n}{2}-dn/2$.
Let $a=2r/n$, and suppose $a\gg\log n$.
Let $F$ be distributed as $\pi_r$.
Then there is an absolute constant $C$ such that, with probability at least
$1-n^{-100}$, simultaneously for every
$e\in E(F)$, 
\[
 |x_e(F)-q|
 \le C\sqrt{\frac{\log n}{a}}\,q(1-q),
\]
where $q=d/(d+a)$.
\end{lemma}

\begin{proof}
Take $\eta=C_0\sqrt{\log n/a}=o(1)$ as in
\Cref{lem:pseudorandom-properties}.
Call $\mathcal T\in\mathcal F_d(F)$ good if both conclusions of that
lemma hold for $(\mathcal T,F)$, and let $b(F)$ be the fraction of
tuples that are not good.
By \Cref{lem:pseudorandom-properties}, we have that $\mathbb{E}[ b(F)]\le n^{-200}$, and thus by Markov's inequality, we get
\begin{align*}
    \mathbb P\bigl(b(F)>n^{-100}\bigr)\le n^{-100}.
\end{align*}
Fix $F$ with $b(F)\le n^{-100}$. The degree condition depends only on
$F$, so it holds whenever $b(F)<1$. Moreover, every tuple satisfies $d_{F\setminus K(\mathcal T)}(v)=d_F(v)-d$ and
\[
 \Delta(F \setminus K(\mathcal T))\le(1+\eta)a.
\]
Fix $e=uv\in E(F)$ and, for convenience, write $x=x_e(F)$,
$Z=|\mathcal F_d(F)|$ and $b=b(F)\le n^{-100}$.
Let $k=\lceil20\log n\rceil$.
Consider the bipartite multigraph $\mathcal B_e$ with parts
\[
 \mathcal M_e
 =\{\mathcal T\in\mathcal F_d(F):e\in K(\mathcal T)\},
 \qquad
 \bar{\mathcal M_{e}}
 =\{\mathcal T\in\mathcal F_d(F):e\notin K(\mathcal T)\}.
\]
Each oriented switching that flips one colour along a simple alternating
$2k$-cycle containing $e$, using edges unused by the other colours (i.e. the other perfect matchings in $\mathcal{T}$),
defines an edge between the original and resulting tuples.
We now double-count the number of edges
$E=e(\mathcal B_e)$, by looking at the edges leaving each side. Note that
\[
 |\mathcal M_e|=xZ,\qquad |\bar{\mathcal M_{e}}|=(1-x)Z.
\]
For a good tuple $\mathcal T=(M_1, \dots, M_d)\in\mathcal M_e$, let $M_i$ be the
matching containing $e=uv$. Every edge of the multigraph $\mathcal B_e$
incident to $\mathcal T$ is obtained by choosing an oriented simple
$M_i$-alternating cycle of length $2k$ containing $e$ and flipping
$M_i$ along this cycle.
If the cycle is oriented starting from $u$, it is counted by
$X_k^{(i)}(u,u;\mathcal T,F)$, and if it is oriented starting from $v$, it
is counted by $X_k^{(i)}(v,v;\mathcal T,F)$. Therefore, by \Cref{lem:pseudorandom-properties}, we have
\[
\deg_{\mathcal B_e}(\mathcal T)
=
X_k^{(i)}(u,u;\mathcal T,F)
+
X_k^{(i)}(v,v;\mathcal T,F)
=
2(1\pm\eta)\frac{a^k}{n}.
\]
For any tuple, including a bad one, we get via the maximum degree bound at most $2((1+\eta)a)^k$ such switchings.
Since there are at most $bZ$ bad tuples, in total we get
\begin{align}
 \label{eq:left-lower}
 (|\mathcal M_e|-bZ)\frac{2(1-\eta)a^k}{n} \le E \le |\mathcal M_e|\frac{2(1+\eta)a^k}{n}
       +2bZ((1+\eta)a)^k.
\end{align}
From the other side, given a good tuple $\mathcal T=(M_1, \dots, M_d) \in\overline{\mathcal M_e}$, we may insert $e$ into
any of the $d$ colours.
Fix $M_i$ the perfect matching of colour $i$. 
Every edge of the multigraph $\mathcal B_e$
incident to $\mathcal T$ and corresponding to inserting $e$ into
$M_i$ is obtained by choosing an oriented simple $M_i$-alternating
cycle of length $2k$ containing $e$ and flipping $M_i$ along this cycle.
If the cycle is oriented so that $e$ is used from $u$ to $v$,
then after removing $e$ and the matching edge incident to $u$, the
remaining part is counted by a path of length $2(k-1)$ from $v$ to the
partner of $u$ in $M_i$. The opposite orientation is counted similarly
by a path of length $2(k-1)$ from $u$ to the partner of $v$ in $M_i$.
Therefore, by \Cref{lem:pseudorandom-properties}, each colour contributes $2(1\pm\eta)\frac{a^{k-1}}{n}$ switchings, and therefore 
\[
\deg_{\mathcal B_e}(\mathcal T)
=
2d(1\pm\eta)\frac{a^{k-1}}{n}.
\]
Moreover, for an arbitrary tuple, this degree is at most $2d((1+\eta)a)^{k-1}$.
Consequently,
\begin{align}
 \label{eq:right-lower}
 (|\bar{\mathcal M_{e}}|-bZ)
       \frac{2d(1-\eta)a^{k-1}}{n} \le E&\le |\bar{\mathcal M_{e}}|
       \frac{2d(1+\eta)a^{k-1}}{n}
       +2dbZ((1+\eta)a)^{k-1}.
\end{align}

Since $k=\lceil20\log n\rceil$ and $\eta=o(1)$, we have
$n(1+\eta)^k\le n^2$ for sufficiently large $n$.
Combining \eqref{eq:left-lower} with \eqref{eq:right-lower}, we get
\[
 \begin{aligned}
  a(1-\eta)(x-b)
  &\le d\bigl((1+\eta)(1-x)+n^2b\bigr),\\
  d(1-\eta)(1-x-b)
  &\le a\bigl((1+\eta)x+n^2b\bigr).
 \end{aligned}
\]
Recall that $q=d/(d+a)$. Rearranging these two inequalities yields
\[
 \frac{q(1-\eta)-n^2b}{1-\eta(2q-1)}
 \le x\le
 \frac{q(1+\eta)+n^2b}{1+\eta(2q-1)}.
\]
Subtracting $q$ and using $\eta\le1/2$, we conclude that
\[
 |x_e(F)-q|
 \le \frac{2\eta q(1-q)+n^2b}{1-\eta}
 \le 4\eta q(1-q)+2n^{-98} = O(\eta q(1-q)).
\]
Since $e$ was arbitrary, the desired estimate holds simultaneously
for every edge of $F$.
\end{proof}

We next show that close to uniform estimates on $x_e(F)$ turn into a coupling with essentially the desired properties.
The following result is based on the coupling approach developed by Gao,
Isaev and McKay \cite{gao2022sandwiching}.

\begin{lemma}
\label{lem:coupling}
Let $N=\binom n2$, $m=dn/2$, $L=N-m$, and $q_r=m/(m+r)$.
Suppose $1\le R<L/2$, and, for every $r>R$, there exist deterministic $0\le\delta_r\le1/2$ and $\beta_r\ge0$ such that 
\[
 \mathbb P_{F\sim\pi_r}\left(
   |x_e(F)-q_r|\le\delta_rq_r(1-q_r)
   \text{ for every }e\in E(F)
 \right)\ge1-\beta_r.
\]
Let $V=\sum_{r=R+1}^{L}\delta_rq_r$.
Then there is a coupling of $F_R\sim\pi_R$ and
$H_*\sim G(n,m+2R)$ such that
\[
 \mathbb P(F_R\nsubseteq H_*)
 \le\sum_{r=R+1}^{L}\beta_r+e^{2V-R}.
\]
\end{lemma}

\begin{proof}
Recall that $\pi_r$ is obtained by taking a uniform ordered $d$-tuple of edge-disjoint perfect matchings and adding a uniform $r$-subset of its complement. Thus, if $F$ has $m+r$ edges, \[ \pi_r(F)=\frac{Z(F)}{Z(K_n)\binom Lr}, \] where $Z(F)$ is the number of ordered $d$-tuples of edge-disjoint perfect matchings contained in $F$. 
We start by defining a deletion process that maps $\pi_r$ to $\pi_{r-1}$. 
Given $F$, delete $e\in E(F)$ with probability 
\[ P_F(e)=\frac{1-x_e(F)}r =\frac{Z(F-e)}{rZ(F)}. \] 
This is a probability distribution since $\sum_{e\in E(F)}Z(F-e)=rZ(F)$, and it is easy to see that this deletion sends $\pi_r$ to $\pi_{r-1}$. Indeed, for any $G$ with $m+r-1$ edges, 
\begin{align*} 
\mathbb P(F\setminus e=G) =\sum_{e\notin E(G)}\pi_r(G+e)P_{G+e}(e) =\frac{(L-r+1)Z(G)} {rZ(K_n)\binom Lr} =\frac{Z(G)}{Z(K_n)\binom L{r-1}} =\pi_{r-1}(G). 
\end{align*}
Call $F$ good if it satisfies $|x_e(F)-q_r|\le\delta_rq_r(1-q_r)$ for every $e\in E(F)$.
Writing $U_F$ for the uniform distribution on $E(F)$ and
$\varepsilon_r=\delta_rq_r$, a good $F$ satisfies
\[
 P_F(e) \geq \frac{1-q_r-\delta_r q_r(1-q_r)}{r} =(1-\varepsilon_r) U_F(e).
\]
Thus, for some probability distribution $Q_F$,
\begin{equation}\label{eq:deletion-mixture}
 P_F=(1-\varepsilon_r)U_F+\varepsilon_rQ_F.
\end{equation}

Start $F=H=K_n$, and independently sample
$B_r\sim\operatorname{Bernoulli}(\varepsilon_r)$ for every $r>R$.
As long as the current graph $F$ is good (we say we are in a good state), we couple one deletion from
$F$ with the uniform deletion process on $H$ as follows:
\begin{itemize}
 \item If $B_r=0$, uniformly delete present edges from $H$ until the
 first edge belonging to $F$ is deleted, and delete that same edge
 from $F$.
 \item If $B_r=1$, delete an edge of $F$ with law $Q_F$, and leave
 $H$ unchanged.
\end{itemize}
If a bad $F$ is encountered, abandon the coupling and continue
$F$ separately with its exact deletion law. 
In any case, we stop once $F$ has $m+R$ edges, and let $F_R$ and $\tilde{H}$ be the obtained graphs. 
Note that by \eqref{eq:deletion-mixture}, after $r$ steps of this process, $F$ has distribution $\pi_{L-r}$, and that for a good state, both updates preserve $F\subseteq H$.
Thus, to prove the desired result, it now suffices to show that the probability of encountering a bad $F$ through the entire process is low, and that $\tilde{H}$ has the required number of edges with high enough probability.
For the first part, the probability of encountering
a bad state before level $R$ is clearly at most $\sum_{r>R}\beta_r$ by union bound.
For the second part, let $B=\sum_{r>R}B_r$. By independence we have
\[
 \mathbb E [e^B]
 =\prod_{r>R}\bigl(1+(e-1)\varepsilon_r\bigr)
 \le e^{2V},
\]
and thus by Markov's inequality
\begin{align*}
    \mathbb P(B>R)\le e^{2V-R}.
\end{align*}
Note that $e(\tilde{H}) - e(F_R)  \leq B$.
Indeed, the number of edges in $H\setminus F$ increases
by at most one when $B_r=1$, and does not increase when $B_r=0$.
Let $H_*$ be the $(m+2R)$-edge level of the completed $H$-process.
Its deletion ordering is uniform, so $H_*\sim G(n,m+2R)$.
If no bad graph is encountered and $B\le R$, then $E(\tilde{H}) - E(F_R)  \leq B \leq R$, and therefore
\[
 F_R\subseteq \tilde{H} \subseteq H_*.
\]
Consequently,
\[
 \mathbb P(F_R\nsubseteq H_*)
 \le \sum_{r>R}\beta_r+e^{2V-R},
\]
as required.
\end{proof}

We are now ready to prove \Cref{thm:Kim-Vu}.
\begin{proof}[Proof of \Cref{thm:Kim-Vu}]
Let $N=\binom n2$ and $m=\frac{dn}{2}$. Let $A$ and $R$ be such that $R = An/2$ and $A= C \sqrt{d \log n}$ for a sufficiently large $C$. Note that $A \gg \log n$, and therefore \Cref{lem:edge-marginals} applies to every $r>R$ with $\beta_r=n^{-100}$ and 
\[
 \delta_r=C_0 \sqrt{\frac{n\log n}{2r}},
\]
for some constant $C_0$ that does not depend on $C$.
Moreover, recalling that $q_r=m/(m+r)$, we have
\begin{align*}
 \sum_{r=R+1}^{L}\delta_rq_r = \sum_{r=R+1}^{L} C_0\sqrt{\frac{n\log n}{2}} \frac{m}{mr^{1/2}+r^{3/2}} &\le C_0\sqrt{\frac{n\log n}{2}}
 \left(\sum_{r=1}^{m}r^{-1/2}
       +m\sum_{r>m}r^{-3/2}\right) \\
       &\le 2C_0n\sqrt{d\log n} \\
       &\le R/4,
\end{align*}
for $C$ large enough.
Since $R<L/2$ for large $n$, where $L=N-m$, and $R = \omega(1)$, by \Cref{lem:coupling} there exists a coupling
of $G_\nu\sim\nu_d$ and $H_*\sim G(n,m+2R)$ such that
$\mathbb P(G_\nu\subseteq H_*)=1-o(1)$.

Let $p_+=\frac{d}{n-1}+\frac{10 C n \sqrt{d \log n}}{\binom{n}{2}}$, and independently sample
$B\sim\operatorname{Bin}(N,p_+)$. Since
\[
 \mathbb E[B] \geq m+3R,
 \qquad \operatorname{Var}(B)=O(m)=o(R^2),
\]
we have $B\ge m+2R$ with probability $1-o(1)$.
The $B$-edge level $G_+$ of the same uniform ordering as $H_*$
has law $G(n,p_+)$ and contains $H_*$ on this event.
Hence $\mathbb P(G_\nu\subseteq G_+)=1-o(1)$.

Finally, by \Cref{thm:Coupling-Decomposition-PerfectMatchings}, there exists a coupling of $G_\nu$
and $G_d\sim\mu_d$ such that $\mathbb{P}(G_d = G_\nu )=1-o(1)$. 
In total, we get a coupling of $G_d \sim\mu_d$ and $G_+ \sim G(n,p_+)$ such that
\[
 \mathbb P(G_d\nsubseteq G_+) \leq \mathbb P(G_d \neq G_\nu) + \mathbb P(G_\nu \not\subseteq G_+) 
 =o(1).
\]
This finishes the proof.
\end{proof}

\bibliographystyle{abbrv}
\bibliography{bib}

@article{hollom2025monotonicity,
  title={Monotonicity and decompositions of random regular graphs},
  author={Hollom, Lawrence and Lichev, Lyuben and Mond, Adva and Portier, Julien and Wang, Yiting},
  journal={The Annals of Applied Probability},
  volume={36},
  number={5},
  pages={4716--4738},
  year={2026},
  publisher={Institute of Mathematical Statistics}
}

@article{IMcKSZ23,
    author={Isaev, M. and McKay, B. D. and Southwell, A. and Zhukovskii, M.},
     TITLE = {Sprinkling with random regular graphs},
    JOURNAL = {Electronic Journal of Probability},
    VOLUME = {30},
      YEAR = {2025},
     PAGES = {1--20},
}

@article{BMW2026,
  title={Asymptotic number of edge-colored regular graphs},
  author={Borinsky, Michael and Meroni, Chiara and Wiesmann, Maximilian},
  journal={arXiv preprint arXiv:2601.18994},
  year={2026}
}

@article{behague2025proof,
  title={A proof of the Kim-Vu sandwich conjecture},
  author={Behague, Natalie and Il'kovi{\v{c}}, Daniel and Montgomery, Richard},
  journal={arXiv preprint arXiv:2510.20765},
  year={2025}
}

@article{Gao23,
  title={The number of perfect matchings, and the nesting properties, of random regular graphs},
  author={Gao, P.},
  journal={Random Structures and Algorithms},
  volume={62},
  number={4},
  pages={935--955},
  year={2023},
  publisher={Wiley Online Library}
}

@article{gao2022sandwiching,
  title={Sandwiching dense random regular graphs between binomial random graphs: P. Gao et al.},
  author={Gao, Pu and Isaev, Mikhail and McKay, Brendan D},
  journal={Probability Theory and Related Fields},
  volume={184},
  number={1},
  pages={115--158},
  year={2022},
  publisher={Springer}
}

@article{gao2025sandwiching,
  title={Sandwiching between random regular graphs and Erd$\backslash$H $\{$o$\}$ sR$\backslash$'enyi graphs: configuration model and unions of perfect matchings},
  author={Gao, Pu and Isaev, Mikhail and Perez-Gimenez, Xavier},
  journal={arXiv preprint arXiv:2510.21472},
  year={2025}
}

@article{hollom2026approximate,
  title={Approximate Itai--Zehavi conjecture for random graphs},
  author={Hollom, Lawrence and Lichev, Lyuben and Mond, Adva and Portier, Julien and Wang, Yiting},
  journal={Random Structures \& Algorithms},
  volume={69},
  number={2},
  pages={e70095},
  year={2026},
  publisher={Wiley Online Library}
}

@article{kim2004sandwiching,
  title={Sandwiching random graphs: universality between random graph models},
  author={Kim, Jeong Han and Vu, Van H},
  journal={Advances in Mathematics},
  volume={188},
  number={2},
  pages={444--469},
  year={2004},
  publisher={Elsevier}
}

@article{dudek2017embedding,
  title={Embedding the Erd{\H{o}}s--R{\'e}nyi hypergraph into the random regular hypergraph and Hamiltonicity},
  author={Dudek, Andrzej and Frieze, Alan and Ruci{\'n}ski, Andrzej and {\v{S}}ileikis, Matas},
  journal={Journal of Combinatorial Theory, Series B},
  volume={122},
  pages={719--740},
  year={2017},
  publisher={Elsevier}
}

@article{klimovsova2023sandwiching,
  title={Sandwiching biregular random graphs},
  author={Klimo{\v{s}}ov{\'a}, Tereza and Reiher, Christian and Ruci{\'n}ski, Andrzej and {\v{S}}ileikis, Matas},
  journal={Combinatorics, Probability and Computing},
  volume={32},
  number={1},
  pages={1--44},
  year={2023},
  publisher={Cambridge University Press}
}

@article{gao2011kim,
  title={Kim--vu’s sandwich conjecture is true for $d \geq (log n)^4$},
  author={Gao, Pu and Isaev, Mikhail and McKay, Brendan},
  journal={arXiv preprint arXiv:2011.09449}
}

@article{mckay1985asymptotics,
  title={Asymptotics for symmetric 0-1 matrices with prescribed row sums},
  author={McKay, Brendan D},
  journal={Ars Combin},
  volume={19},
  pages={15--25},
  year={1985}
}

@article{janson1990poisson,
  title={Poisson approximation for large deviations},
  author={Janson, Svante},
  journal={Random Structures \& Algorithms},
  volume={1},
  number={2},
  pages={221--229},
  year={1990},
  publisher={Wiley Online Library}
}

@article {MRRW97,
    AUTHOR = {Molloy, M. S. O. and Robalewska, H. and Robinson, R. W. and
              Wormald, N. C.},
     TITLE = {{$1$}-factorizations of random regular graphs},
   JOURNAL = {Random Structures and Algorithms},
    VOLUME = {10},
      YEAR = {1997},
    NUMBER = {3},
     PAGES = {305--321},
      ISSN = {1042-9832},
MRREVIEWER = {Mark R. Jerrum},
}

@article{janson1995random,
  title={Random regular graphs: asymptotic distributions and contiguity},
  author={Janson, Svante},
  journal={Combinatorics, Probability and Computing},
  volume={4},
  number={4},
  pages={369--405},
  year={1995},
  publisher={Cambridge University Press}
}

\end{document}